\documentclass{amsart}
\usepackage{amssymb, amsmath, amsfonts, amscd}

\input xypic

\theoremstyle{plain}
\newtheorem{theorem}{Theorem}[section]
\newtheorem{corollary}[theorem]{Corollary}
\newtheorem{lemma}[theorem]{Lemma}
\newtheorem{proposition}[theorem]{Proposition}

\newtheorem{example}[theorem]{Example}

\theoremstyle{definition}
\newtheorem{definition}[theorem]{Definition}

\theoremstyle{remark}

\numberwithin{equation}{theorem}

\newcommand{\F}{\mathcal{F}}

\newcommand{\I}{\mathcal{I}}
\newcommand{\J}{\mathcal{J}}

\newcommand{\E}{\mathcal{E}}

\newcommand{\Hom}{\operatorname{Hom}} 
\renewcommand{\O}{\mathcal{O} }

\renewcommand{\P}{\mathbf{P} }
\renewcommand{\Pr}{\mathcal{J} }

\newcommand{\Exan}{\operatorname{Exan}}
\newcommand{\Diff}{\operatorname{Diff}}

\newcommand{\Ext}{\operatorname{Ext} }

\newcommand{\Der}{\operatorname{Der} }

\newcommand{\Real}{\mathbf{R}}
\newcommand{\R}{\mathcal{R}}
\newcommand{\KO}{\operatorname{KO}}
\newcommand{\K}{\operatorname{K}}
\newcommand{\spec}{\operatorname{Spec}}

\newcommand{\p}{\mathbf{P}}

\begin{document}

\title{On jets, extensions  and characteristic classes II}

\author{Helge Maakestad}
%\address{Institut Fourier, Grenoble}

\email{\text{h\_maakestad@hotmail.com} }


\date{Spring 2009}

\begin{abstract} 
In this paper we define the generalized Atiyah classes $c_\Pr(\E)$ and $c_{\O_X}(\E)$ of a
quasi coherent sheaf $\E$ with respect to a pair $(\I,d)$ where $\I$
is a left and right $\O_X$-module and $d$ a derivation. 
We relate this class to the structure of left and right module on the first
order jet bundle $\Pr^1_\I(\E)$. In the main Theorem of the paper 
we show $c_{\O_X}(\E)=0$ if and only if there is an isomorphism
$\Pr^1_\I(\E)^{left} \cong \Pr^1_\I(\E)^{right}$ as $\O_X$-modules. 
We also give explicit examples where $c_{\O_X}(\E)\neq 0$ using jet bundles of line
bundles on the projective line hence the classes $c_\Pr(\E)$ and
$c_{\O_X}(\E)$ are non trivial. 
The classes we introduce generalize the classical Atiyah class.
\end{abstract}

\maketitle

\tableofcontents

\section{Introduction} 

The aim of this paper is to introduce and study generalized first
order jet bundles and generalized Atiyah classes for quasi coherent
sheaves relative to an arbitrary morphism $\pi:X\rightarrow S$ of
schemes. We define the generalized first order jet bundle $\Pr^1_\I(\E)$ of $\E$
and the generalized Atiyah sequence
\begin{align}
&\label{at} 0\rightarrow \I\otimes_{\O_X}\E \rightarrow \Pr^1_\I(\E) \rightarrow \E \rightarrow 0.
\end{align}
The generalized Atiyah sequence is an exact sequence of quasi coherent
sheaves of left and right $\O_X$-modules and left and right $\Pr^1_\I$-modules. The sheaf $\Pr^1_\I$ is a
sheaf of associative rings on $X$. It is an extension of $\O_X$ with
a left and right quasi coherent  $\O_X$-module $\I$ of square zero. 
The main result of the paper is that the sequence \ref{at} is left
split as sequence of $\O_X$-modules if and only if the left and right structure on $\Pr^1_\I(\E)$
are $\O_X$-isomorphic (see Theorem \ref{main}).

We give a general definition of the first order jet bundle
$\Pr^1_\I(\E)$ of a
quasi coherent sheaf $\E$ using a derivation $d$ and a left
and right $\O_X$-module $\I$. We define generalized Atiyah classes
$c_{\Pr}(\E)$ and $c_{\O_X}(\E)$ of $\E$ and relate these classes to the left and right
$\O_X$ structure on $\Pr^1_\I(\E)$. The generalized Atiyah class
$c_{\O_X}(\E)$ measures when there is
an isomorphism
\[ \Pr^1_\I(\E)^{left}\cong \Pr^1_\I(\E)^{right} \]
as $\O_X$-modules (see Theorem \ref{main}).
There is always an isomorphism
\[ \Pr^1_\I(\E)^{left}\cong \Pr^1_\I(\E)^{right} \]
as sheaves of abelian groups. Hence the class $c_{\O_X}(\E)$ 
measures when the sheaf of abelian groups $\Pr^1_\I(\E)$ may be given two non
isomorphic structures as $\O_X$-module.
When $\I=\Omega^1_{X/S}$ and $d:\O_X\rightarrow \Omega^1_{X/S}$ is the
universal derivation it follows the characteristic class
$c_{\O_X}(\E)$ is the classical Atiyah class as defined in \cite{atiyah}.

We prove in Example \ref{line} there is no isomorphism
\[ \Pr^1_{\p^1}(\O(d))^{left}\neq \Pr^1_{\p^1}(\O(d))^{right} \]
of $\O_{\p^1}$-modules on $\p^1$. Hence the class
$c_{\O_{\p^1}}(\O(d))$ is non-zero for the
sheaf $\O(d)$ on the projective line $\p^1$ over a field $K$ of
characteristic zero. If $char(K)$ divides $d$ it follows
$c_{\O_{\p^1}}(\O(d))$ is zero hence in this case there is an isomorphism
\[ \Pr^1_{\p^1}(\O(d))^{left}\cong \Pr^1_{\p^1}(\O(d))^{right} \]
of $\O_{\p^1}$-modules. It follows the generalized Atiyah classes we
introduce are non trivial.

In a previous paper on jet bundles (see \cite{maa2})
we found many examples of jets of line bundles on the projective line where
the structure as left locally free $\O_{\p^1}$-module was nonisomorphic to the 
structure as right locally free $\O_{\p^1}$-module. In this paper we interpret this
phenomenon in terms of the generalized Atiyah class.

In the first section of the paper we do all constructions in the local
case. In the second part of the paper we do all constructions in the
global case. In the final section of the paper we discuss Atiyah
sequences and Atiyah classes relative to arbitrary morphisms of
schemes and arbitrary infinitesimal extensions of $\O_X$ by a quasi
coherent left and right $\O_X$-module. We get a
definition of Atiyah classes $c_\Pr(\E)$ and $c_{\O_X}(\E)$ for en 
arbitrary quasi coherent sheaf $\E$ of $\O_X$-modules. Here $\J$ is 
an infinitesimal extension of $\O_X$ with respect to a quasi coherent left and right $\O_X$-module $\I$.
The construction generalize the classical case (see Example \ref{genrclassical}).
In the final section of the paper we construct explicit examples. We
also construct (see Example \ref{complex}) for any smooth projective scheme $X$ over the complex
numbers a characteristic class $c_{\I}(\E)$ in $\K_0(X)$ where
$\K_0(X)$ is the Grothendieck group of locally free finite rank
sheaves on $X$. There is a canonical map of Grothendieck groups
\[ \gamma:\K_0(X)\rightarrow \KO(X(\Real)) \]
where $X(\Real)$ is the underlying real smooth manifold of $X$ and
$\KO(X(\Real))$ is the Grothendieck group of the category of real
finite rank smooth vector bundles on $X(\Real)$.
The class $c_{\I}(\E)$ lies in the group $ker(\gamma)$ for any
locally free finite rank sheaf $\E$ and any left and right
$\O_X$-module $\I$.

\section{Generalized Atiyah classes for modules over rings}

Let in the following $\phi:A\rightarrow B$ be a unital morphism of
commutative rings. Let $E$ be a $B$-module. Let $I$ be a left and
right $B$-module with $a(xb)=(ax)b$ for all $a,b\in B$ and $x\in I$
and let $d\in \Der_A(B,I)$ be a derivation. This means $d(\phi(a))=0$ for
all $a\in A$ and $d(ab)=ad(b)+d(a)b$ for all $a,b\in B$. In the
following we will write $d(a)$ instead of $d(\phi(a))$.

Let $\Pr^1_I=I\oplus B$ with the following multiplication:
\[ (x,a)(y,b)=(xb+ay,ab).\]
It follows $\Pr^1_I$ is an associative ring with multiplicative unit
$\mathbf{1}=(0,1)$. Define a map $d_I:B\rightarrow \Pr^1_I$ by
$d_I(b)=(d(b),0)$. Define maps $s,t:B\rightarrow \Pr^1_I$ by
\[ t(b)=(0,b) \]
and
\[ s(b)=(d(b),b).\]
Define two left actions of $B$ on $\Pr^1_I$ as follows:
\begin{align}
&\label{taction}(x,a)_tb=(x,a)t(b)=(xb,ab)\text{ and } b_t(x,a)=t(b)(x,a)=(bx,ba)
\end{align}
and
\begin{align}
&\label{saction}(x,a)_sb=(x,a)s(b)=(xb+ad(b),ab)\\
&\label{saction1}b_s(x,a)=s(b)(x,a)=(d(b)a+bx,ba).
\end{align}
The ring $\Pr^1_I$ acts canonically on $I\subseteq \Pr^1_I$ as follows:
\[ (x,0)(y,b)=(xb,0)\text{ and }(y,b)(x,0)=(bx,0).\]

In the following when we write $\Der_A^t(B,\Pr^1_I)$ we mean 
$A$-linear derivations with $B$-structure on $\Pr^1_I$ induced by $t$
from \ref{taction}.
\begin{proposition} The following holds:
\begin{align}
&\label{ass1}\Pr^1_I\text{ is a square zero extension of $B$ by
  $I$.}\\
&\label{ass2}\text{\ref{taction}-\ref{saction1} define $\Pr^1_I$ as
  left and right $B$-module}.\\
&\label{ass3}\text{The map $d_I$ is a derivation: $d_I\in
  \Der_A^t(B,\Pr^1_I)$.}\\
&\label{ass4}\text{The maps $s,t$ are ring homomorphisms and
  $s-t=d_I$.}
\end{align}
\end{proposition}
\begin{proof} The proof is left to the reader.
\end{proof}

Define the following map:
\[ d_\Pr:\Pr^1_I\rightarrow I\]
by
\[ d_\Pr(x,a)=x+d(a).\]

When we write $\Der_A^t(\Pr^1_I,I)$ we mean derivations linear over
$A$ with respect to the left action of $A$ on $\Pr^1_I$ induced by $t$
from \ref{taction}.
\begin{lemma} The following holds: $d_\Pr\in \Der_A^t(\Pr^1_I,I)$.
\end{lemma}
\begin{proof} The proof is left to the reader.
\end{proof}

Define the following abelian group:
\[ \Pr^1_I(E)=\Pr^1_I\otimes_B E\cong I\otimes_B E\oplus E.\]

Define the following left and right action of $\Pr^1_I$ on
$\Pr^1_I(E)$:
\begin{align} 
&\label{a1}(x,a)(z\otimes e,f)=(x\otimes f+az\otimes e+d(a)\otimes f, af) 
\end{align}
and
\begin{align}
&\label{a2} (z\otimes e,f)(x,a)=(z\otimes ea, fa).
\end{align}

Note: We may write
\[ (x,a)(z\otimes e,f)=((az)\otimes e+d_\Pr(x,a)\otimes f, af).\]

\begin{definition} Let $\Pr^1_I(E)$ with the actions \ref{a1} and
\ref{a2} be the \emph{first order $I$-jet module} of $E$ with respect to $d$.
\end{definition}

There is an exact sequence of abelian groups
\begin{align}
&\label{atiyah}0\rightarrow I\otimes_B E\rightarrow \Pr^1_I(E)\rightarrow E \rightarrow 0
\end{align}
There is a left and right action of $\Pr^1_I$ on $E$ defined as
follows:
\begin{align}
&\label{e1}(x,a)e=ae
\end{align}
\begin{align}
&\label{e2}e(x,a)=ea
\end{align}

\begin{proposition} \label{modules} The actions \ref{a1} and \ref{a2} make
  $\Pr^1_I(E)$ into a left and right $\Pr^1_I$-module. The exact
  sequence \ref{atiyah} is an exact sequence of left and right $\Pr^1_I$-modules.
\end{proposition}
\begin{proof} The proof is left to the reader.
\end{proof}

\begin{corollary} The sequence \ref{atiyah} is an exact sequence of
  left and right $B$-modules. It is split exact as right $B$-modules.
\end{corollary}
\begin{proof} Use the ring homomorphism $t:B\rightarrow \Pr^1_I$
  defined by $t(b)=(0,b)$. It follows the maps in the sequence
  \ref{atiyah} are $B$-linear. The Corollary now follows from
  Proposition \ref{modules}.
\end{proof}

\begin{definition} Let the sequence \ref{atiyah} be the
  \emph{Atiyah-Karoubi sequence} of $E$ with respect to the pair $(I,d_\Pr)$.
\end{definition}

Note: The exact sequence \ref{atiyah} was first defined in \cite{atiyah} in the case when
$\I=\Omega_{X/\mathbf{C}}$,  $d$ is the universal derivation and $X$
is a complex manifold.

Recall the following definition: Given two left $B$-modules $E,F$ we define
\[ \Diff^0_A(E,F)=\Hom_B(E,F)\]
and
\[ \Diff^k_A(E,F)=\{ \partial \in \Hom_A(E,F): [\partial,a]\in
\Diff^{k-1}_A(E,F)\text{ for all }a\in B\}.\]

There is a natural map
\[ d_E:E\rightarrow \Pr^1_I(E) \]
defined by
\[ d_E(e)=(0,1)\otimes e\]
called the \emph{universal differential operator of $E$} with respect
to the pair $(I,d)$. Here $(0,1)=\mathbf{1}\in \Pr^1_I=I\oplus B$ is
the multiplicative unit.

\begin{lemma} \label{universal}  It follows $d_E\in \Diff_A^1(E,\Pr^1_I(E))$.
\end{lemma}
\begin{proof}
Recall the following: For elements $(x,a)\in \Pr^1_I$ and $b\in B$ it
follows
\[ (x,a)_tb=(x,a)(0,b)=(xb ,ab) \]
and
\[ b_s(x,a)=(d(b),b)(x,a)=(bx+d(b)a,ba) .\]

It follows $d_E\in \Diff^1(E,\Pr^1_I(E))$ if and only if $[d_E,a]\in
\Hom_B(E,\Pr^1_I(E))$ for all $a\in B$. We get
\[ [d_E,a](be)=(d_Ea-ad_E)(be)=d_E(abe)-ad_E(be)=(0,1)\otimes abe
-a(0,1)\otimes be=\]
\[ (0,ab)\otimes e-a(0,1)\otimes be=(0,ab)\otimes e-(d(a)b,ab)\otimes
e=\]
\[(-d(a)b,0)\otimes e.\]
By definition
\[ a(0,1)=(d(a),a)=(d(a),0)+(0,a)=(d(a),0)+(0,1)a.\]
It follows
\[ (-d(a),0)=(0,1)a-a(0,1).\]
We get
\[ [d_E,a](be)=(-d(a)b,0)\otimes e=b(-d(a),0)\otimes
e)=b((0,1)a\otimes e-a(0,1)\otimes e)=\]
\[b((0,1)\otimes ae-a(0,1)\otimes e)=b(d_E(ae)-ad_E(e))=b[d_E,a](e).\]
Hence $[d_E,a]\in \Hom_B(E,\Pr^1_I(E))=\Diff^0(E,\Pr^1_I(E))$ for all
$a\in B$.
It follows $d_E\in \Diff^1_A(E,\Pr^1_I(E)) $ and the claim of the
Lemma follows.
\end{proof}

Recall there is a derivation $d_\Pr\in \Der_A(\Pr^1_I,I)$ defined by
$d_\Pr(x,a)=x+d(a)$. It follows $d_\Pr|_B=d_I$.
Let $F$ be a left $\Pr^1_I$-module. An $A$-linear map
\[ \nabla:F\rightarrow I\otimes_B F \]
satisfying
\[ \nabla((x,a)f)=(x,a)\nabla(f)+d_\Pr(x,a)\otimes f \]
is an $(I,d_\Pr)$-connection on $F$. 
Assume $E$ is a $B$-module. An $A$-linear map
\[ \nabla:E\rightarrow I\otimes_B E \]
satisfying 
\[ \nabla(ae)=a\nabla(e)+d(a)\otimes e \]
is an $(I,d)$-connection.

\begin{lemma} \label{connection} Assume $\nabla$ is an $(I,d_\Pr)$-connection. It follows 
$\nabla$ is an $(I,d)$-connection.
\end{lemma}
\begin{proof} Let $b\in B$ and $f\in F$. We get
\[ \nabla(bf)=\nabla((0,b)f)=(0,b)\nabla(f)+d_\Pr(0,b)\otimes f=\]
\[ b\nabla(f)+d(b)\otimes f \]
and the Lemma is proved.
\end{proof}

The sequence $\ref{atiyah}$ is an exact sequence of left
$\Pr^1_I$-modules. We get an extension class
\[ c_\Pr(E)\in \Ext^1_{\Pr^1_I}(E,I\otimes_B E).\]
When we restrict to $B$ we get an extension class
\[ c_B(E)\in \Ext^1_B(E,I\otimes_B E) .\]

\begin{definition} Let $c_\Pr(E)$ be the \emph{Atiyah class of $E$ with
    respect to $\Pr^1_I$}. Let $c_B(E)$ be the \emph{Atiyah class of
    $E$ with respect to $B$}.
\end{definition}

\begin{proposition} \label{split} The following holds:
\begin{align}
&\label{at1}\text{\ref{atiyah} is split as sequence of right
  $\Pr^1_I$-modules.}\\
&\label{at2}\text{\ref{atiyah} is left split as $\Pr^1_I$-modules iff $E$ has an $(I,d_\Pr)$-connection.}\\
&\label{at3}\text{\ref{atiyah} is left split as $B$-modules iff $E$ has an $(I,d)$-connection.}
\end{align}
\end{proposition}
\begin{proof} We prove claim \ref{at1}: Define the following map:
\[ s:E\rightarrow \Pr^1_I(E) \]
by
\[ s(e)=(0,e).\]
It follows 
\[s(e(x,a))=s(ea)=(0,ea)=(0,e)(x,a)=s(e)(x,a) \] 
and claim \ref{at1} is proved. We prove claim \ref{at2}:
Assume $s(e)=(\nabla(e),e)$ is a left $\Pr^1_I$-linear section. We get
\[ s((x,a)e)=(\nabla((x,a)e),(x,a)e)=(\nabla((x,a)e), ae)=\]
\[ (x,a)(\nabla(e),e)=(x\otimes e+a\nabla(e)+d(a)\otimes e,
ae)=((x,a)\nabla(e)+(x+d(a))\otimes e, ae).\]
It follows $\nabla$ satisfies
\[ \nabla((x,a)e)=(x,a)\nabla(e)+d_\Pr(x,a)\otimes e \]
and $\nabla$ is a $(I,d_\Pr)$-connection. Claim \ref{at2} follows.
Claim \ref{at3} follows in a similar way and the Proposition is proved.
\end{proof}

\begin{example} The classical Atiyah class. \end{example}

When $I=\Omega_{B/A}=\Omega$ and $d:B\rightarrow \Omega$ is the
universal derivation it follows the class $c_B(E)$ is the classical
Atiyah  class as defined in \cite{atiyah}. From Proposition
\ref{split} it follows $c_{B}(E)=0$ if and only if $E$ has a
connection
\[ \nabla:E\rightarrow \Omega_{B/A}\otimes_B E \]

\begin{corollary} If $c_\Pr(E)=0$ it follows $c_B(E)=0$.
\end{corollary}
\begin{proof} If $c_\Pr(E)=0$ it follows from Proposition \ref{split},
  claim \ref{at2} $E$ has an $(I,d_\Pr)$-connection. It follows from
  Lemma \ref{connection} $E$ has a $(I,d)$-connection. From this and
  Proposition \ref{split}, claim \ref{at3} the Lemma follows.
\end{proof}

Let $\Pr^1_I(E)^{left}$ denote the abelian group $\Pr^1_I(E)$ with its
left $B$-module structure. 
Let $\Pr^1_I(E)^{right}$ denote the abelian group $\Pr^1_I(E)$ with its
right $B$-module structure. We say the module $I$ is
\emph{abelianized} if the following holds:
\[ ax =xa \]
for all $a\in B$ and $x\in I$. Define the following product on $I$:
\[ a*x=xa\]
for all $a\in B$ and $x\in I$. It follows we have defined an
$B$-module structure on $I$, denoted $I^{star}$ with the property
there is an isomorphism
\[ I^{right}\cong  I^{star} \]
of $A$-modules. When we form the tensor product
\[ I\otimes_B E \]
we use the right structure on $I$ and left structure on $E$. Since $B$
is commutative it follows $E$ has a canonical right $B$-module
structure.
The abelian group $I\otimes_B E$ has a left and right structure as
$B$-module, denoted
\[ I\otimes_B E^{left} \]
and
\[ I\otimes_B E^{right}.  \]
We form a new product on $I\otimes_B E$ as follows:
\[ a*(x\otimes e)=x\otimes (ea) \]
for any $a\in B$ and $x\otimes e\in I\otimes_B E$. We get a left
$B$-module denoted $I\otimes_B E^{star}$. There is an isomorphism
\[ I\otimes_B E^{star}\cong I\otimes_B E^{right} \]
of $B$-modules.

\begin{lemma} \label{abelianized} Assume $I$ is abelianized. It follows there is an
  isomorphism
\[ I\otimes_B E^{left} \cong I\otimes_B E^{right} \]
of $B$-modules.
\end{lemma}
\begin{proof} Define the following map:
\[ \phi:I\otimes_B E^{left} \rightarrow I\otimes_B E^{star} \]
by
\[ \phi(x\otimes e)=x\otimes e.\]
We get
\[\phi(a(x\otimes e))=\phi((ax)\otimes e)=(ax)\otimes e=(xa)\otimes
e=\]
\[ x\otimes (ae) =x\otimes (ea)=a*(x\otimes e)=a*\phi(x\otimes e).\]
It follows
\[I\otimes_B E^{left}\cong I\otimes_B E^{star}\cong I\otimes_B
E^{right}\]
and the Lemma follows.
\end{proof}

\begin{proposition} Assume $I$ is abelianized. \label{equiv} The following holds:
\begin{align}
&\label{split1}c_B(E)=0\text{ iff $\Pr^1_I(E)^{left}\cong \Pr^1_I(E)^{right}$ as $B$-modules.}\\
&\label{split2}\Pr^1_I(E)^{left}\cong \Pr^1_I(E)^{right}\text{ as
  abelian groups.}\\
&\label{split3}\Pr^1_I(E)^{left} \cong \Pr^1_I(E)^{right}\text{ as $B$-modules iff
  $E$ has an $(I,d)$-connection}
\end{align}
\end{proposition}
\begin{proof} It follows $c_B(E)=0$ if and only if the sequence
\[ 0\rightarrow I\otimes_B E \rightarrow \Pr^1_I(E) \rightarrow E
\rightarrow 0 \]
is split as sequence of left  $B$-modules. From Lemma
\ref{abelianized} it follows
\[ \Pr^1_I(E)^{left}\cong I\otimes_B E\oplus E^{left} \cong \]
\[ I\otimes_A E\oplus E^{right} \cong \Pr^1_I(E)^{right} \]
since by Proposition \ref{split}, claim \ref{at1} the sequence is
always split as right $B$-modules. Claim \ref{split1} follows.
Claim \ref{split2} is obvious. Claim \ref{split3}: From \ref{split1}
we get
\[ \Pr^1_I(E)^{left}\cong \Pr^1_I(E)^{right}\]
as $B$-modules iff $c_B(E)=0$. From Proposition \ref{split}, claim
\ref{at3} it follows $c_B(E)=0$ iff $E$ has an
$(I,d)$-connection. Claim \ref{split3} is proved and the Proposition follows.
\end{proof}

Hence the characteristic class $c_B(E)$ measures when the abelian
group $\Pr^1_I(E)$ may be equipped with two non isomorphic structures
as $B$-module.

\section{Generalized Atiyah classes for quasi coherent sheaves} 

In this section we generalize the results in the previous section to
the case where we consider an arbitrary morphism $\pi:X\rightarrow S$
of schemes and an arbitrary quasi coherent $\O_X$-module $\E$. We
define the generalized Atiyah class $c_\Pr(\E)$ and $c_{\O_X}(\E)$
using derivations and sheaves of associative rings and prove various
properties of this construction. We end the section with a discussion
of explicit examples. We give examples of jet bundles
$\Pr^1_\Omega(\O(d))$ where $\O(d)=\O(1)^{\otimes d}$ and $\O(1)$ is
the tautological quotient bundle on $\P^1_K$. We prove
$c_{\O_{\p^1}}(\O(d))\neq 0$ when $char(K)=0$ and
$c_{\O_{\p^1}}(\O(d))=0$ when $char(K)$ divides $d$. Hence the
characteristic classes $c_{\Pr}$ and $c_{\O_X}$ are non trivial.

Let in the following $\pi:X\rightarrow S$ be an arbitrary morphism of
schemes and let $\I$ be an arbitrary quasi coherent
left and right $\O_X$-module. Let
\[ d\in \Der_{\pi^{-1}(\O_S)}(\O_X,\I) \]
be a derivation. We may form the quasi coherent sheaf
$\Pr^1_\I=\I\oplus \O_X$. Let $V\subseteq U$ be open subsets of $X$
and let $(x,a),(y,b)\in \Pr^1_\I(U)$ be two elements. Define
\[ (x,a)(y,b)=(xb+ay,ab).\]
It follows $\Pr^1_\I(U)$ is an associative ring with multiplicative
unit $\mathbf{1}=(0,1)$ and the restriction morphism
\[ \Pr^1_\I(U)\rightarrow \Pr^1_\I(V) \]
is a morphism of unital rings. There is a natural embedding
\[ i:\I\rightarrow \Pr^1_\I \]
defined as follows: 
\[ i(U):\I(U)\rightarrow \Pr^1_\I(U) \]
\[ i(U)(x)=(x,0).\]
It follows $\I\subseteq \Pr^1_\I$ is a sheaf of ideals in $\Pr^1_\I$
with $\I^2=0$.
We get an exact sequence
\begin{align}\label{squarezero}
 0\rightarrow \I\rightarrow   \Pr^1_\I  \rightarrow \O_X \rightarrow 0
\end{align}
of sheaves of abelian groups.
The sequence \ref{squarezero} is an extension of $\O_X$ by a quasi
coherent sheaf of two sided ideals of square zero.

We may
form the quasi coherent sheaf $\Pr^1_\I(\E)=\I\otimes_{\O_X} \E\oplus
\E$. We may define a left and right $\O_X$-structure
and a left $\Pr^1_\I$-module structure on $\Pr^1_\I(\E)$ as follows:
Let $(x,a)\in \Pr^1_\I(U)$ and $(z\otimes e,f)\in
\Pr^1_\I(\E)(U)$. Let
\[ (x,a)(z\otimes e, f)=(x\otimes f+az\otimes e+d(a)\otimes f, af).\]
and
\[(z\otimes e,f)(x,a)=(z\otimes ea,fa).\]
One checks for any open sets $V\subseteq U$ it follows the restriction morphism
\[\rho_{UV}: \Pr^1_\I(\E)(U)\rightarrow \Pr^1_\I(\E)(V) \]
satisfies
\[ (x,a)(z\otimes e,f)|_V=(x,a)|_V(z\otimes e,f)|_V .\]
For open sets $W\subseteq V\subseteq U$ it follows 
\[ \rho_{VW}\circ \rho_{UV}=\rho_{UW}.\]
Similar formulas holds for the right structure as $\Pr^1_\I$-module.
It follows $\Pr^1_\I(\E)$ becomes a sheaf of left and right
$\Pr^1_\I$-modules  and left and right $\O_X$-modules.
We get an exact sequence
\begin{align}\label{global}
0\rightarrow \I\otimes_\O \E\rightarrow \Pr^1_\I(\E)\rightarrow
\E\rightarrow 0
\end{align}
of sheaves left and right $\O_X$-modules and left and right $\Pr^1_\I$-modules.
These assertions follow immediately from the local situation since all
sheaves involved are quasi coherent. 

\begin{definition} Let the sequence \ref{global} be the \emph{Atiyah-Karoubi
    sequence of $\E$ with respect to $(\I,d)$}.
\end{definition}

Define the following $\pi^{-1}(\O_S)$-linear map:
\[ d_\E:\E \rightarrow \Pr^1_\I(\E) \]
by
\[ d_\E(U)(e)=(0,1)\otimes e \in \Pr^1_\I(\E)(U).\]
Here $U\subseteq X$ is an open subset.

\begin{lemma} The following holds:
\[ d_\E\in \Diff^1_{\pi^{-1}(\O_S)}(\E,\Pr^1_\I(\E)).\]
\end{lemma}
\begin{proof} The Lemma follows from Lemma \ref{universal} since the sheaf $\E$ is quasi coherent.
\end{proof}
The morphism $d_\E$ is the \emph{universal differential operator for
  $\E$ with respect to $(\I,d)$}.

\begin{example}The first order jet bundle. \end{example}

Assume $\I=\Omega^1_{X/S}$ and $d$ the universal derivation. We get a
map
\[ d_\E:\E\rightarrow \Pr^1_{X/S}(\E) \]
defined by
\[ d_\E(U)(e)=1\otimes e\in \Pr^1_{X/S}(\E)(U).\]
This map is the classical differential operator $d_\E\in
\Diff^1_{\pi^{-1}(\O_S)}(\E,\Pr^1_{X/S}(\E))$ for $\E$.

In the following we view the sequence \ref{global} as an exact sequence
of sheaves of left $\Pr^1_\I$ and $\O_X$-modules.
We get a characteristic class
\[ c_\Pr(\E)\in \Ext^1_{\Pr^1_\I}(\E,\I\otimes_{\O_X} \E) \]
and
\[ c_{\O_X}(\E)\in \Ext^1_{\O_X}(\E,\I\otimes_{\O_X} \E) \]

\begin{definition}
Let $c_\Pr(\E)$ be the \emph{Atiyah class of $\E$ with respect to $\Pr^1_\I$}.
Let $c_{\O_X}(\E)$ be the \emph{Atiyah class of $\E$ with respect to $\O_X$}.
\end{definition}

There is a derivation $d_\Pr\in \Der_{\pi^{-1}(\O_S)}(\Pr^1_\I,\I)$
defined as follows: Let $(x,a)\in \Pr^1_\I(U)$. Define
\[ d_\Pr(x,a)=x+d(a).\]

Assume $\E$ is a quasi coherent $\O_X$-module.
We say an $\pi^{-1}(\O_S)$-linear map
\[ \nabla:\E\rightarrow \I\otimes_{\O_X} E \]
is an \emph{($\I$,d)-connection} if for all local sections
$a\in \O(U)$ and $e\in \E(U)$ on an open set $U\subseteq X$ the
following holds:
\[ \nabla(ae)=a\nabla(e)+d(a)\otimes e .\]

Assume $\F$ is a left $\Pr^1_\I$-module which is quasi coherent as
left $\O_X$-module.
We say an $\pi^{-1}(\O_S)$-linear map
\[ \nabla:\F\rightarrow \I\otimes_{\O_X} \F \]
is an $(\I,d_\Pr)$-connection if the following holds for
$(x,a)\in \Pr^1_\I(U)$ and $f\in \F(U)$:
\[ \nabla((x,a)f)=(x,a)\nabla(f)+d_\Pr(x,a)\otimes f.\]

\begin{proposition}\label{equality}  The following holds:
\begin{align}
&\label{c1}c_\Pr(\E)=0\text{ iff $\E$ has an $(\I,d_\Pr)$-connection.}\\
&\label{c2}c_{\O_X}(\E)=0\text{ iff $\E$ has an $(\I,d)$-connection.}\\
&\label{c3}\text{If $c_\Pr(\E)=0$ it follows $c_{\O_X}(\E)=0$.}
\end{align}
\end{proposition}
\begin{proof} The proof is left to the reader.
\end{proof}

We say the sheaf $\I$ is \emph{abelianized} if for all local sections
$a$ of $\O_X$ and $x$ of $\I$ the following holds
\[ ax=xa.\]
As in Lemma \ref{abelianized} We get an isomorphism
\[ I\otimes_{\O_X}\E \oplus \E^{left} \cong I\otimes_{\O_X}\E \oplus
\E^{right} \]
of $\O_X$-modules.

\begin{theorem}Assume $\I$ is abelianized. \label{main} The following holds:
\begin{align}
&\label{m1}\text{The sequence \ref{global} is split as sequence of
  right $\O_X$-modules} \\
&\label{m11}\text{\ref{global} is split as left $\Pr^1_\I$-module iff
  $\E$ has an $(\I,d_\Pr)$-connection.}\\
&\label{m4}\Pr^1_\I(\E)^{left}\cong \Pr^1_\I(\E)^{right}\text{ as
  sheaves of abelian groups.}\\
&\label{m3}c_{\O_X}(\E)=0\text{ iff }\Pr^1_\I(\E)^{left}\cong
\Pr^1_\I(\E)^{right}\text{ as $\O_X$-modules.}\\
&\label{m2}\text{\ref{global} is split as left $\O_X$-modules iff $\E$
  has an $(\I,d)$-connection} 
\end{align}
\end{theorem}
\begin{proof} We prove \ref{m1}: Let $U\subseteq X$ be an open subset
  and define the morphism
\[ s(U):\E(U)\rightarrow \Pr^1_\I(\E)(U) \]
by
\[ s(U)(e)=(0,e).\]
Let $X=(x\otimes a,b)\in \Pr^1_\I(U)$.
It follows
\[ s(U)(eX)=(0,eb)=(0,e)X=s(U)(e)X\]
hence $s$ is left $\Pr^1_\I$-linear. The map $s$ splits \ref{global}
as sequence of left $\Pr^1_\I$-modules and claim \ref{m1} follows.
We prove \ref{m11}: By definition \ref{global} is split as left
$\Pr^1_\I$-modules iff $c_\Pr(\E)=0$. By proposition this is iff $\E$
has an $(\I,d_\Pr)$-connection. Claim \ref{m11} is proved.
We prove \ref{m2}: The sequence \ref{global} is split as left
$\O_X$-modules iff $c_{\O_X}(\E)=0$. By Proposition \ref{equality}
this is iff $\E$ has an $(\I,d)$-connection. Claim \ref{m2} is proved.
We prove \ref{m3}: By definition it follows 
\[ \Pr^1_\I(\E)^{right}\cong \I\otimes_{\O_X}\E\oplus \E^{right} \]
as right $\O_X$-modules. Sequence \ref{global} is left split as
$\O_X$-modules iff 
\[ \Pr^1_\I(\E)^{left}\cong \I\otimes_{\O_X}\E\oplus \E^{left} .\]
It follows $c_{\O_X}(\E)=0$ iff there is an isomorphism
\[ \Pr^1_\I(\E)^{left}\cong \I\otimes_{\O_X}\E\oplus \E^{left}\cong
\I\otimes_{\O_X}\E\oplus \E^{right} \cong \Pr^1_\I(\E)^{right}.\]

Claim \ref{m3} follows.
Claim \ref{m4} is obvious and the Theorem follows.
\end{proof}

Hence the characteristic class $c_{\O_X}(\E)$ measures when the sheaf of
abelian groups $\Pr^1_\I(\E)$ is equipped with two non isomorphic
structures as $\O_X$-module.

\begin{example} The classical case: $\I=\Omega^1_{X/S}$. \end{example}

Assume in the following Proposition $\pi:X\rightarrow S$ is a
separated morphism.
Let $\Delta:X\rightarrow X\times_S X$ be the diagonal embedding. It
follows $\Delta(X)\subseteq X\times_S X$ is a closed subscheme. Let
$\J\subseteq \O_{X\times_S X}$ be the ideal sheaf of $\Delta(X)$ and
let $\O_{\Delta^l}=\O_{X\times_S X}/\J^{l+1}$ be the \emph{l'th
infinitesimal neighborhood of the diagonal}. Let
$p,q:X\times_S X\rightarrow X$ be the canonical projection maps.

\begin{definition} Let $\Pr^l_{X/S}(\E)=p_*(\O_{\Delta^{l+1}}\otimes
      q^*\E)$ be the \emph{$l$'th order jet bundle of $\E$}.
\end{definition}

Assume $\pi$ is given by a homomorphism $\phi:A\rightarrow B$ of
commutative rings. Let $X=\spec(B)$ and $S=\spec(A)$. Assume $\E$ is
the sheaffification of a $B$-module $E$.
It follows
$\Pr^1_{X/S}(\E)$ is the sheaf associated to $P^l_{B/A}(E)=B\otimes_A
B/J^{l+1}\otimes_B E$ where $J\subseteq B\otimes_A B$ is the kernel of
the multiplication map.

\begin{proposition}  \label{recover} Assume $\I=\Omega_{X/S}$ and $d:\O_X\rightarrow
  \Omega_{X/S}$ is the universal derivation. It follows
  $\Pr^1_\I(\E)\cong \Pr^1_{X/S}(\E)$ is the first order jet bundle of
  $\E$. The generalized Atiyah sequence \ref{global} becomes the
  classical Atiyah sequence.
\end{proposition}
\begin{proof} Assume $U=\spec(B)\subseteq X$ is an open affine
  subscheme mapping to an open affine subscheme $V=\spec(A)\subseteq
  S$. Let $\E(U)=E$ where $E$ is a $B$-module and let
  $\Omega_{X/S}|(U)=\Omega$.
Let $m:B\otimes_A B\rightarrow B$ be the multiplication map and
  let $s:B\rightarrow P^1_{B/A}$ be defined by $s(b)=1\otimes b$. It
  follows there is an isomorphism
\[ \phi:P^1_{B/A}\otimes_B E\rightarrow \Omega\otimes_B E\oplus E \]
defined by
\[ \phi(x\otimes e)=((x-sm(x))\otimes e, m(x)e) .\]
One checks we get a commutative diagram of exact sequences
\[
\diagram 0 \rto & \Omega\otimes_B E \rto \dto^= & P^1_{B/A}(E) \rto
\dto^\cong & E \dto^= \rto & 0 \\
  0 \rto & \Omega\otimes_B E \rto & \Pr^1_\Omega(E) \rto & E \rto & 0
\enddiagram \]
where the middle vertical arrow is the isomorphism $\phi$. Since the
map $\phi$ is intrinsically defined it glues to give an isomorphism
$\Pr^1_{X/S}(\E)\cong \Pr^1_\Omega(\E)$. The Proposition follows.
\end{proof}

Note: If $I=\Omega_{B/A}$ and $d:B\rightarrow \Omega_{B/A}$ is the universal
derivation the
construction of $\Pr^1_{\Omega_{B/A}}(E)$ using $d$ is due to
Karoubi (see \cite{karoubi}).
In this case we get the sequence
\[ 0\rightarrow \Omega_{B/A}\otimes_B E\rightarrow
P^1_{B/A}(E)\rightarrow E \rightarrow 0 \]
where $P^1_{B/A}(E)$ is the \emph{first order jet module} of $E$. It
is also called the \emph{first order module of principal parts} of
$E$.

Note: Assume $\E$ is a finite rank locally free $\O_X$-module.
If $X$ is a smooth scheme of finite type over $\mathbf{C}$ - the
complex numbers, $\I=\Omega^1_{X}$ is the module of differentials and
$d$ the universal derivation we get the classical Atiyah sequence
\[ 0\rightarrow \Omega^1_X\otimes \E \rightarrow
\Pr^1_{\Omega^1_X}(\E)\rightarrow \E \rightarrow 0.\]
The class 
\[ c_{\O_X}(\E)\in \Ext^1_{\O_X}(\E,\Omega^1_X\otimes \E) \]
is the classical Atiyah class as defined in \cite{atiyah}.
It follows $c_{\O_X}(\E)=0$ iff $\E$ has a connection
\[ \nabla:\E\rightarrow \Omega^1_X\otimes \E.\]

\begin{example} \label{line} Sheaves of jets on the projective
  line.\end{example}

Let $K$ be a field of characteristic zero and consider $\p^1_K$. Let
$\O(d)$ be $\O(1)^{\otimes d}$ where $\O(1)$ is the tautological
quotient bundle on $\p^1_K$. Let $\Omega=\Omega_{\p^1_K}^1$ be the sheaf of
differentials on $\p^1_K$.
We get from \ref{recover} an exact sequence of
$\O_{\p^1}$-modules
\[ 0\rightarrow \Omega \otimes \O(d)\rightarrow
\Pr^1_\Omega(\O(d))\rightarrow \O(d)\rightarrow 0 .\]
It follows $\Pr^1_\Omega(\O(d))\cong \Pr^1_{\p^1}(\O(d))$ is the first
order jet bundle of $\O(d)$ on $\p^1_K$.
There is a left and right $\O_{\p^1}$-module structure on
$\Pr^1_{\p^1_K}(\O(d))$ and by \cite{maa2} isomorphisms
\[ \Pr^1_{\p^1_K}(\O(d))^{left}\cong \O(d-1)\oplus \O(d-1) \]
and
\[ \Pr^1_{\p^1_K}(\O(d))^{right}\cong \O(d)\oplus \O(d-2) \]
of $\O_{\p^1}$-modules. It follows from Theorem \ref{main}, claim
\ref{m3} $c_{\O_{\p^1}}(\O(d))\neq 0$. There is an ismorphism
\[ \Pr^1_{\p^1_K}(\O(d))^{left}\cong \Pr^1_{\p^1_K}(\O(d))^{right} \]
as sheaves of abelian groups.

By Proposition \ref{equality} it follows $c_\Pr(\O(d))\neq 0$ if
$char(K)=0$
hence the classes $c_{\O_{\p^1}}(\O(d))$ and  $c_\Pr(\O(d))$ are non trivial.

Again by \cite{maa2} if $char(K)$ divides $d$ it follows there is an
isomorphism
\[ \Pr^1_{\p^1_K}(\O(d))^{left}\cong \Pr^1_{\p^1_K}(\O(d))^{right} \]
of $\O_{\p^1}$-modules. Hence in this case $c_{\O_{\p^1}}(\O(d))=0$. 
The sheaf $\Pr^1_{\p^1}(\O(d))$ and the class $c_{\O_{\p^1}}(\O(d))$ is defined
over
$\p^1_{\mathbf{Z}}$. When we pass to $\mathbf{Q}$ the class $c_{\O_{\p^1}}(\O(d))$
is non zero. When we pass to $\mathbf{F}_p=\mathbf{Z}/p\mathbf{Z}$ and
when $p$ divides $d$ the class $c_{\O_{\p^1}}(\O(d))$ becomes zero.

\section{Jets and infinitesimal extensions of sheaves}

In this section we define and study jet bundles $\Pr^1_\I(\E)$ of $\O_X$-modules $\E$
relative to any morphism $\pi:X\rightarrow S$ of schemes. We define in
Definition \ref{genrseq} the generalized Atiyah sequence of $\E$ with
respect to a square zero extension $\J$ and a derivation $d\in
\Der_\R(\O_X,\I)$.  Here $\J$ is any sheaf of associative unital
algebras on $X$ that are square zero extensions of $\O_X$ with respect
to a quasi coherent left and right $\O_X$-module $\I$.
We use the generalized Atiyah sequence
\ref{genrseq} in Definition \ref{genratiyah} to define the generalized
Atiyah classes $c_\J(\E)$ and $c_{\O_X}(\E)$ for $\E$. This construction
generalize the classical construction (see Example \ref{genrclassical}).

Let in the following $\pi:X\rightarrow S$ be a morphism of schemes and
let $\R=\pi^{-1}(\O_S)$. 
It follows $\O_X$ is a sheaf of $\R$-algebras. This means for any open subset
$U\subseteq X$ the ring $\O_X(U)$ is an $\R(U)$-algebra.
Let $\E$ be a quasi coherent sheaf of left $\O_X$-modules. 
It follows $\E$ is canonically a sheaf of quasi coherent right $\O_X$-modules.
Let 
\[ 0\rightarrow \I\rightarrow \J \rightarrow^m \O_X \rightarrow 0\]
be an infinitesimal exstension of $\O_X$ by a quasi coherent left and right
$\O_X$-module $\I$. This means $\J$ is a sheaf of associative unital
rings on $X$  and $\I\subseteq \J$ is a sheaf of quasi coherent two sided ideals with
$\I^2=0$. Let $d\in \Der_\R(\O_X,\I)$ be a derivation over $\R$. This
means for any open set $U$ $d(U)\in \Der_{\R(U)}(\O_X(U),\I(U))$ is a
derivation over $\R(U)$. The derivation $d$ induce a derivation
$p=d\circ m$ in $\Der_\R(\J,\I)$.
Make the following definition
\[ \Pr^1_\I(\E)=\I\otimes_{\O_X}\E\oplus \E .\]
Let $a\in \J(U)$ and $Z=(z\otimes e,f)\in \Pr^1_\I(\E)(U)$ and define
\[ aZ=a(z\otimes e,f)=((az)\otimes e+p(a)\otimes f,af) \]
and
\[ Za=(z\otimes e,f)a=(z\otimes (ea),fa).\]
this is well defined since the sheaf $\O_X$ is a sheaf of commutative
rings hence $\E$ is canonically a sheaf of right $\O_X$-modules.
It follows for any open subset $V\subseteq U$ there is an equality
\[ a(z\otimes e,f)|_V=a|_V(z\otimes e,f)|_V.\]
If $\rho_{UV}$ is the restriction map from 
\[ \rho_{UV}:\Pr^1_\I(\E)(U)\rightarrow \Pr^1_\I(\E)(V) \]
it follows for any open sets $W\subseteq V \subseteq U$ we get
\[\rho_{VW}\circ \rho_{UV}=\rho_{UW}.\]
Let $W=(w\otimes g,h)\in \Pr^1_\I(\E)(U)$ and $b\in \J(U)$.
\begin{proposition} The sheaf $\Pr^1_\I(\E)$ is a sheaf of left and right
  $\J$-modules. The sheaf $\I\otimes_{\O_X}\E$ is a sheaf of
left and right  $\J$-submodules of $\Pr^1_\I(\E)$. The sheaf $\E$ is a
sheaf of left and right $\J$-modules.
\end{proposition}
\begin{proof} The proof is left to the reader.
\end{proof}

Note: $\Pr^1_\I(\E)$ is not a lifting of $\E$ to $\J$ in the sense of
deformation theory since
$\I\Pr^1_\I(\E)=0$.
It follows
\[ \O_X \otimes_\J \Pr^1_\I(\E)\cong (\J/\I)\otimes_\J
\Pr^1_\I(\E)\cong \]
\[ \Pr^1_\I(\E)/\I\Pr^1_\I(\E)\cong \Pr^1_\I(\E)\neq \E .\]
By definition a lifting $\E_\J$ of $\E$ to $\J$ is required to satisfy
\[ \O_X \otimes_\J \E_\J\cong \E.\]

We get a natural exact sequence of sheaves of abelian groups
\begin{align}
&\label{ringed}0\rightarrow \I\otimes_{\O_X}\E \rightarrow^i \Pr^1_\I(\E) \rightarrow^j \E \rightarrow 0.
\end{align}

\begin{definition} \label{genrseq} Let the sequence \ref{ringed} be the
  \emph{Atiyah-Karoubi sequence of $\E$} with respect to the pair $(\I,d)$.
\end{definition}

\begin{corollary} The sequence \ref{ringed} is an exact sequence of
  left and right $\J$-modules and left and right $\O_X$-modules.
\end{corollary}
\begin{proof} One checks the natural maps $i,j$ are left and right
  $\J$-linear, left and right $\O_X$-linear and the Corollary follows.
\end{proof}

View the sequence \ref{ringed} as a sequence of left $\J$ and
$\O_X$-modules.
We get two characteristic classes:
\[ c_\J(\E)\in \Ext^1_\J(\E,\I\otimes_{\O_X}\E) \]
and
\[ c_{\O_X}(\E)\in \Ext^1_{\O_X}(\E,\I\otimes_{\O_X}\E) .\]

\begin{definition} \label{genratiyah} The class $c_\J(\E)$ is the \emph{Atiyah class of
    $\E$} with respect to $\J$. The class $c_{\O_X}(\E)$ is the
  \emph{Atiyah class of $\E$} with respect to $\O_X$.
\end{definition}

\begin{example} \label{genrclassical} The classical Atiyah class.\end{example}

 The class $c_{\O_X}(\E)$ generalize the classical Atiyah class: When
$X$ is a separated scheme over a fixed base scheme $S$ and
$\I=\Omega^1_{X/S}$ is the sheaf of differentials, $d$ is the
universal derivation, $\J=\I\times \O_X$ is the trivial square zero extension
and $\R=\pi^{-1}(\O_S)$ it follows $c_{\O_X}(\E)$
is the classical Atiyah class as defined in \cite{atiyah}.

\section{Appendix: Explicit examples}

In this section we give some explicit examples to illustrate the
constructions made in the previous sections. Let in the following $K$
be an arbitrary field.

\begin{example} Connections on projective modules.\end{example}

Let $A$ be a commutative ring with unit and let $P$ be a fintely
generated projective $A$-module. Let 
\[ 0\rightarrow K \rightarrow A^n \rightarrow^p P \rightarrow 0\]
be an exact sequence where $A^n$ is the free $A$-module of rank $n$.
Since $P$ is projective the map $q$ has a section $s$ with $q\circ s =id$.

Let $F=A^n=A\{e_1,..,e_n\}$.
Let $p_i:F\rightarrow A$ be defined by $p_i=e_i^*$. Let $q_i=p_i\circ
s$. 
We get elements
\[ e_1,..,e_n\in P\]
and
\[ q_1,..,q_n\in P^* \]
satisfying the following:
\[ \sum_{i=1}^n q_i(a)e_i=a \]
for any $a\in P$.
Define the following map:
\[ \nabla(e)=\sum_{i=1}^n d(q_i(e))\otimes e_i .\]

\begin{lemma} \label{connection} The map $\nabla$ is a connection on
  $P$.
\end{lemma}
\begin{proof} The map $\nabla$ is by definition well-defined. We check
  it is a connection.
\[ \nabla(ae)=\sum_i d(q_i(ae))\otimes e_i = \sum_i d(aq_i(e))\otimes e_i =\]
\[ \sum_i (ad(q_i(e))+d(a)q_i(e))\otimes e_i= a\sum_id(q_i(e))\otimes e_i +d(a)\otimes \sum_iq_i(e)e_i=\]
\[ a\nabla(e)+d(a)\otimes e \]
and the map $\nabla$ is a connection. The Lemma follows.
\end{proof}

The universal derivation $d:A\rightarrow \Omega^1$ gives rise to the
Atiyah-Kaorubi sequence
\[ 0\rightarrow \Omega^1\otimes P \rightarrow
\Pr^1_{\Omega^1}(P)\rightarrow P \rightarrow 0.\]
The abelian group $\Pr^1(P)=\Pr^1_{\Omega^1}(P)$ has a canonical
left $A$-structure and a canonical right $A$-structure. It has by the
previous section a
non-trivial left $A$-structure defined as follows:
\[ a(x\otimes e, f)=((ax)\otimes e+d(a)\otimes f, af) \]
with $a\in A$ and $ (x\otimes e,f)\in \Pr^1(P)$.
By Lemma \ref{connection} and Theorem \ref{main} there is an
isomorphism
\[ \Pr^1(P)^{left} \cong \Pr^1(P)^{right} \]
of $A$-modules.

\begin{example} Generalized Atiyah classes on the projective line\end{example}

In this example we consider generalized Atiyah-sequences of locally
free sheaves on the projective line. We prove that in the case of an
invertible sheaf the only non-trivial Atiyah-class arise in the
classical case. The aim of the constructions made is to construct new
examples of locally free sheaves of left and right modules on the
projective line where the splitting type as left module differs from the splitting
type as right module.

Let $V=K\{e_0,e_1\}$ and $V^*=K\{x_0,x_1\}$ and $\p=\p(V^*)$ the
projective line over $K$. Let $\O(m)=\O(1)^{\otimes m}$ be the $m$'th
tensor power of the tautological line bundle on $\p$. There is an
isomorphism
\[ \Omega^1_{\p}\cong \O(-2) \]
where $\Omega^1_{\p}$ is the cotangent bundle on $\p$. We get an exact
sequence
\[ 0\rightarrow \Omega^1_{\p}\otimes \O(m) \rightarrow \J_\p(\O(m))
\rightarrow \O(m) \rightarrow 0\]
called the \emph{classical Atiyah sequence} of $\O(m)$.
Assume
\[d\in  \Der_K(\O_\p ,\O(m)) \]
is a derivation and consider the generalized Atiyah-sequence
\begin{align}
&\label{ati} 0\rightarrow \O(m)\otimes \O(l) \rightarrow
\J_{\O(m)}(\O(l))
\rightarrow \O(l) \rightarrow 0
\end{align}
with respect to $(\O(m), d)$. We get a characteristic class
\[ c_{\O_\p}(\O(l))\in \Ext^1(\O(l),\O(m+l)).\]

\begin{proposition} The only non-trivial case is when $m=-2$,
\[ d:\O_\p \rightarrow \Omega^1_\p \]
is the universal derivation and \ref{ati} is the classical Atiyah
sequence.
\end{proposition}
\begin{proof} The proof is left to the reader.
\end{proof}

When we use a derivation $d:\O_\p \rightarrow \I$ where $\I$ is a
higher rank locally free sheaf we may get new examples.

Let $\I=\J_\p\cong \O(-2)\oplus \O$ be the first order jet bundle on
$\p$.
Define the following derivation $d:\O_\p\rightarrow \I$. Let
$U_i=D(x_i)$ and let $t=x_1/x_0$, $s=1/t$ be local coordinates on
$\p$.
Let $\O(U_0)=K[t]$, $\I(U_0)=K[t]\{dt,e\}$ and
\[ d_0:K[t]\rightarrow K[t]\{dt, e\} \]
be defined by
\[ d_0(a(t))=a'(t)dt+t^ia'(t)e\]
with $i=0,1,2$.
Let $\O(U_1)=K[s]$, $\I(U_1)=K[s]\{ds,f\}$ and
\[ d_1:K[s]\rightarrow K[s]\{ds,f\} \]
be defined by
\[ d_1(b(s))=b'(s)ds^-s^{2-i}b'(s)f.\]
One checks $d_0,d_1$ glue to a derivation $d\in \Der(\O_\p, \I)$.
We get for any linebundle $\O(l)$ on $\p$ an Atiyah-Karoubi sequence
\[ 0\rightarrow I\otimes \O(l) \rightarrow \J_\I(\O(l)) \rightarrow
\O(l) \rightarrow 0.  \]
The local structure of $\J_\I(\O(l))$ looks as follows:
\[ \J_\I(\O(l))(U_0)=K[t]\{dt\otimes x_0^l, e\otimes x_0^l\} \]
with left $K[t]$-multiplication given as follows:
\[ a(\omega, bx_0^l)=(a\omega+d_0(a)\otimes bx_0^l, abx_0^l).\]
Moreover
\[ \J_\I(\O(l))(U_1)=K[s]\{ds\otimes x_1^l, f\otimes x_1^l\} \]
with left $K[s]$-multiplication given as follows:
\[ c(\eta, dx_1^l)=(c\eta+d_1(c)\otimes dx_1^l, cdx_1^l).\]
We aim to construct the structure matrix of the locally free sheaf
$\J_\I(\O(l))$ as left $\O_\p$-module and to see how this matrix
depend on the derivation $d$ and the integer $l$. The splitting type
of $\J_\I(\O(l))$ as right $\O_\p$-module is by the previous section
as follows:
\[ \J_\I(\O(l))^{right}\cong \O(l-2)\oplus \O(l)\oplus \O(l).\]
Using Atiyah-Karoubi sequences we hope to give new examples where the
splitting type as left module differs from the one as right module.
Let $\J_\I(\O(l))(U_0)$ have the following basis as left
$K[t]$-module:
\[B: (1,0,0)_0=dt\otimes x_0^l, (0,1,0)_0=e\otimes x_0^l,
(0,0,1)_0=x_0^l .\]
Let $\J_\I(\O(l))(U_1)$ have the following basis as left
$K[s]$-module:
\[C: (1,0,0)_1=ds\otimes x_1^l, (0,1,0)_1=f\otimes x_1^l,
(0,0,1)_1=x_1^l.\]

\begin{theorem} The structure matrix $[L]^C_B$ is as follows:
\[[L]^C_B=
\begin{pmatrix} -t^{l-2} &     0      & -lt^{l-1}      \\
                    0    &     t^l    & -lt^{i+l-1}    \\
                    0    &     0      & t^l
\end{pmatrix}
\]
\end{theorem}
\begin{proof}  The proof is a straight forward calculation.
\end{proof}

We see the structure matrix $[L]^C_B$ depend on the integers
$i,l$. Maybe one will get examples with more than two different
structures of $\O_\p$-module on the locally free sheaf $\O(l-2)\oplus
\O(l)\oplus \O(l)$.

\begin{example} \label{complex} Vector bundles on real and complex
  manifolds. \end{example}

Note:  If $K$ is the field of complex numbers and $X$ a smooth
projective scheme of finite type over $K$, let $X(\mathbf{R})$ denote
the underlying real smooth manifold of $X$.
Let $\KO(X(\mathbf{R}))$ be the Grothendieck group of the category
of real smooth finite rank vectorbundles on $X(\mathbf{R})$.
There is a canonical map
\[ \gamma:\K_0(X)\rightarrow \KO(X(\mathbf{R})) \]
sending the class of a locally free finite rank $\O_{X}$-module $\E$
to the class of it's underlying
finite rank real vectorbundle $\E(\mathbf{R})$. In many cases the
exact sequence
\[ 0\rightarrow \I\otimes \E \rightarrow \Pr^1_{\I}(\E)\rightarrow \E
\rightarrow 0\]
is split as sequence of real smooth vectorbundles on $X(\mathbf{R})$
since $X(\mathbf{R})$ is a real compact manifold and all short exact
sequences of vector bundles split. It follows there is an isomorphism
\[ \Pr^1_{\I}(\E)(\mathbf{R})^{left} \cong
\Pr^1_{\I}(\E)(\mathbf{R})^{right} \]
of real smooth vectorbundles. 
Let 
\[ c_{\I}(\E)=[\Pr^1_{\I}(\E)^{left}]-[\Pr^1_{\I}(\E)^{right}]\in \K_0(X).\]
It follows there is an equality
\[ \gamma(c_{\I}(\E))=0 \]
in $\KO(X(\mathbf{R}))$.

\begin{example} The classical Atiyah sequence.\end{example}

Let $\E=\O(d_1)\oplus \cdots \oplus \O(d_e)$ be a locally free sheaf
on $\P^1$ of rank $e$ over any field $K$.  Define
\[ deg(\E)=d_1+\cdots +d_e.\]

\begin{lemma} \label{projectivebundle} There is an isomorphism
\[ \phi:\K_0(\P^1)\cong \mathbf{Z}\oplus \mathbf{Z} \]
defined by
\[ \phi([\E])=(deg(\E), rk(\E)).\]
\end{lemma}
\begin{proof} By the projective bundle formula there is an isomorphism
\[ \K_0(\P^1)\cong \mathbf{Z}[h]/(1-h)^2 \]
of rings. We get
\[ \phi([\O(-1)]^d)=h^d=(1+h-1)^d=1+d(h-1)=1-dt \]
where $t$ is the class of $1-h$ in $\K_0(\P^1)$. 
We get similarly
\[ \phi([\O(1)]^d)=1-d(h-1)=1+d(1-h)=1+dt.\]
It follows
\[ \phi([\E])=1+d_1t+\cdots +1+d_et=rk(\E)+deg(\E)t \]
and the claim of the Lemma follows.
\end{proof}

Let $\Pr^k(\O(d))$ be the $k$'th sheaf of principal parts of $\O(d)$
on $\P^1$ when $1\leq k \leq d$. By \cite{maa2} it follows
\[ deg(\Pr^k(\O(d))^{left})=deg(\Pr^k(\O(d))^{right}) .\]

\begin{corollary} There is an equality
\[ c_{\Omega}(\O(d))=0\]
in $\K_0(\P^1)$
\end{corollary}
\begin{proof} By Lemma \ref{projectivebundle} it follows
\[ c_{\Omega}(\O(d))=[ \Pr^1(\O(d))^{left}]-[\Pr^1(\O(d))^{right}]=\]
\[ (rk(\Pr^1(\O(d))^{left}),
deg(\Pr^1(\O(d))^{left}))- \]
\[ (rk(\Pr^1(\O(d))^{right}), deg(\Pr^1(\O(d))^{right}))=0.\]
The Corollary follows.
\end{proof}

The triviality of the class $c_\Omega(\O(d))$ is related to the fact
that the first order jet bundle $\Pr^1(\O(d)$ is an extension of two
abelianized $\O_{\P^1}$-modules: $\Omega^1_{\P^1}\otimes \O(d)$ and
$\O(d)$:
the Atiyah sequence
\[ 0\rightarrow \Omega^1_{\P^1}\otimes \O(d)\rightarrow
\Pr^1(\O(d))\rightarrow \O(d) \rightarrow 0 \]
is an exact sequence of left and right $\O_{\P^1}$-modules and there are
isomorphisms
\[ (\Omega^1_{\P^1}\otimes \O(d))^{left} \cong (\Omega^1_{\P^1}\otimes
\O(d))^{right} \]
and
\[ \O(d)^{left} \cong \O(d)^{right} \]
of $\O_{\P^1}$-modules.
Hence
\[
c_\Omega(\O(d))=[\Pr^1(\O(d))^{left}]-[\Pr^1(\O(d))^{right}]=\]
\[ [(\Omega^1_{\P^1}\otimes \O(d))^{left}]+[\O(d)^{left}]-
[(\Omega^1_{\P^1}\otimes \O(d))^{right}]-[\O(d)^{right}] =0\].

This argument holds in the following general case:

\begin{theorem} There is an equality
\[ c_\Omega(\E)=0 \]
in  $\K_0(X)$ where $X/S$ is differentially smooth and $\E$ locally
free of finite rank.
\end{theorem}
\begin{proof} The Atiyah sequence is exact as left and right
  $\O_X$-modules and there are isomorphisms
\[ (\Omega^1_X\otimes \E)^{left}\cong (\Omega^1_X\otimes
\E)^{right} \]
of $\O_X$-modules.
\end{proof}

This example motivates the construction given in this paper: The left and
right $\O$-module $\I$ must be non-abelianized for the class
$c_\I(\E)$ to be non-trivial in $\K_0(X)$. The class $c_\I(\E)$ lies in an
Ext-group which is computable hence it should be easy to check if
there is an isomorphism
\[ \Pr^1_\I(\E)^{left} \cong \Pr^1_\I(\E)^{right} \]
of $\O_X$-modules. 

There is ongoing work where the
characteristic class $c_\I(\E)$ is studied (see \cite{maa4}).

\end{document}